\documentclass[12pt]{article}

\usepackage{graphicx}
\usepackage{nameref}
\usepackage[shortlabels,inline]{enumitem}
\usepackage{empheq}
\usepackage[shortlabels]{enumitem}
\setlist[enumerate]{nosep}
\usepackage[doc]{optional}
\usepackage{xcolor}
\usepackage[colorlinks=true,
            linkcolor=refkey,
            urlcolor=lblue,
            citecolor=red]{hyperref}
\usepackage{float}
\usepackage{soul}
\usepackage{graphicx}
\definecolor{labelkey}{rgb}{0,0.08,0.45}
\definecolor{refkey}{rgb}{0,0.6,0.0}
\definecolor{Brown}{rgb}{0.45,0.0,0.05}
\definecolor{lime}{rgb}{0.00,0.8,0.0}
\definecolor{lblue}{rgb}{0.5,0.5,0.99}

 \usepackage{mathpazo}

\colorlet{hlcyan}{cyan!30}

\usepackage{stmaryrd}
\usepackage{amssymb}
\makeatletter
\def\namedlabel#1#2{\begingroup
   \def\@currentlabel{#2}%
   \label{#1}\endgroup
}
\makeatother

\newcommand{\seppthree}{\setlength{\itemsep}{-3pt}}

\newcommand{\To}{\ensuremath{\rightrightarrows}}

\newcommand{\fenv}[1]%
{\ensuremath{\,\overrightarrow{\operatorname{env}}_{#1}}}
\newcommand{\benv}[1]%
{\ensuremath{\,\overleftarrow{\operatorname{env}}_{#1}}}

\newcommand{\scal}[2]{\left\langle{#1},{#2}  \right\rangle}

\newcommand{\RR}{\ensuremath{\mathbb R}}

\newcommand{\dom}{\ensuremath{\operatorname{dom}}}

\newcommand{\sri}{\ensuremath{\operatorname{sri}}}

\newcommand{\ran}{\ensuremath{{\operatorname{ran}}\,}}
\newcommand{\zer}{\ensuremath{\operatorname{zer}}}

\newcommand{\Id}{\ensuremath{\operatorname{Id}}}

{\begin{list}{}{%
\settowidth{\labelwidth}{\textrm{#1~}}%
\setlength{\leftmargin}{\labelwidth+\labelsep}}}
{\end{list}}
\usepackage{amsthm}
\usepackage[capitalize,nameinlink]{cleveref}
\crefname{equation}{}{equations}
\crefname{chapter}{Appendix}{chapters}
\crefname{item}{}{items}
\crefname{enumi}{}{}
\newtheorem{theorem}{Theorem}[section]
\newtheorem{lemma}[theorem]{Lemma}

\newtheorem{definition}[theorem]{Definition}

\newtheorem{ex}[theorem]{Example}

\newtheorem{remark}[theorem]{Remark}

\providecommand{\norm}[1]{\lVert#1\rVert}

\providecommand{\RR}{\mathbb{R}}

\providecommand{\ran}{\operatorname{ran}}

\providecommand{\dom}{\operatorname{dom}}

\providecommand{\gra}{\operatorname{gra}}
\providecommand{\Id}{\operatorname{{ Id}}}

\providecommand{\To}{\rightrightarrows}

\providecommand{\ran}{\operatorname{ran}}

\providecommand{\Id}{\operatorname{Id}}

\providecommand{\zer}{\operatorname{zer}}

\providecommand{\RR}{\mathbb{R}}

\definecolor{myblue}{rgb}{.8, .8, 1}
  \newcommand*\mybluebox[1]{%
    \colorbox{myblue}{\hspace{1em}#1\hspace{1em}}}

\allowdisplaybreaks 

\begin{document}

%

\title{\textsc{ On Parallel Sums of Set-Valued Operators: Monotonicity, Paramonotonicity, and Decomposition Invariance}}
\author{
Salihah Thabet Alwadani\thanks{
Mathematics, Royal Comission of Jubail and Yanbu,
Yanbu, Saudi Arabia. E-mail:
\texttt{salihah.s.alwadani@gmail.com} and \texttt{wadanis@rcjy.edu.sa}.}
}

\date{November  1 , 2026}
\maketitle

\vskip 8mm

\begin{abstract} 
In this paper, we investigate structural properties of the parallel sum of set-valued operators in a real Hilbert space, with particular emphasis on monotonicity, maximal monotonicity, and paramonotonicity. We examine the decomposition-invariance property of the parallel sum of two paramonotone operators and provide a direct proof based on monotonicity and the equality of admissible decompositions. This approach establishes the required orthogonality relations without invoking the Attouch--Th{\'e}ra Duality duality framework. We also revisit the classical maximal-monotonicity property of the parallel sum and present an alternative graph-based derivation involving an auxiliary shifted operator. Furthermore, we study the parallel sum of a set-valued operator and the normal cone operator associated with a closed linear subspace. Explicit characterizations of its values, domain, and range are presented, and a decomposition-independent representation is obtained under paramonotonicity through the orthogonality of the underlying subspaces. These results provide a unified and self-contained treatment of selected properties of parallel sums and clarify the role of paramonotonicity in their structural behavior.
\end{abstract} 

{\small

}

\maketitle


\section{Preliminaries}

Throughout, we assume that 
\begin{equation}
\text{$ \mathcal{H}$ is
	a real Hilbert space with inner product
	$\scal{\cdot}{\cdot}\colon \mathcal{H}\times \mathcal{H} \to\RR$, }
\end{equation}
and induced norm $\|\cdot\|\colon \mathcal{H} \to \RR\colon x\mapsto \sqrt{\scal{x}{x}}$. Let $A : \mathcal{H} \To \mathcal{H}$ be an arbitrary set valued operator, i.e., $Ax \subseteq \mathcal{H} \  (\forall x \in \mathcal{H})$. The \emph{graph} of $A$ denoted $\gra A$, is defined as $   \gra A= \{  (x, y) \in \mathcal{H} \times \mathcal{H}  \mid y \in Ax  \}. $
Moreover, the \emph{domain} of $A$, denoted by $ \dom A$, is defined as $   \dom A= \{  x \in \mathcal{H} \mid Ax \neq \emptyset \} .$ The \emph{range} of $A$, denoted by $\ran A$, is defined as $    \ran A= \{  y \in \mathcal{H}  \mid \exists x \in \mathcal{H} \ \ \text{such that} \ \ y \in Ax \}. $ The set valued \emph{inverse} $A^{-1}$ of $A$ is defined by $   (y, x)  \in \gra A^{-1} \Leftrightarrow  (x, y) \in \gra A.$ The set of \emph{zeros} of $A$ is written as $  \zer A = A^{-1} (0) = \{  x \in \mathcal{H} \mid 0 \in Ax \}. $ The notation employed in this paper is standard and closely aligns with that in \cite{Alwadanithesis}, \cite{alwadani2021difference}, \cite{alwadani2021attouch}, and \cite{BC2017}.
 
 The key findings are summarized as follows:
\begin{enumerate}
\item \cref{Lem1} establishes a decomposition-invariance property for the parallel sum of two paramonotone operators. Its significance lies in showing that the nonempty intersection of the operator values remains unchanged under different decompositions of a fixed element. This result is closely related to the invariance property established by Bauschke, see \cite[Theorem~5.3]{bauschke2012attouch} in the framework of Attouch–Théra duality. While their proof invokes an orthogonality relation obtained in \cite[Corollary~2.13]{bauschke2012attouch}, our proof derives the required orthogonality relations directly from the monotonicity of the two operators and the equality of the decompositions. Paramonotonicity then yields the equality of the corresponding intersections. 
\item  Our treatment of the parallel sum includes a direct graph-based argument involving the shifted operator $E=(A \square B)-\Id$. This approach provides an alternative derivation of the classical maximal-monotonicity property of the parallel sum under the assumption that $A^{-1} + B^{-1}$ is maximally monotone. Although the maximal-monotonicity result itself is well established, the present argument offers a self-contained formulation based on the shifted operator and the characterization of maximal monotonicity, while maintaining a consistent proof framework throughout the paper, see \cref{Fact01}.
\item \cref{theorem111} investigates the parallel sum of a set-valued operator and the normal cone operator associated with a closed linear subspace. Building on established properties of parallel sums and normal cone operators, we provide explicit characterizations of its values, domain, and range, and examine the role of paramonotonicity in obtaining a decomposition-independent representation. Our approach relies on the orthogonality of the underlying subspaces and provides a direct and self-contained treatment of these properties within the framework of parallel sums.
\end{enumerate}

\begin{definition}\label{def1} \normalfont\cite[Definition~2.18]{Alwadanithesis}
Let $C$ be a nonempty convex subset of $\mathcal{H}$ and $x \in \mathcal{H}$. Then the \emph{normal cone operator} to $C$ at $x$ is 
$$  N_{C} : \mathcal{H} \To \mathcal{H}: x \to 
  \begin{cases}
    \{ y \in \mathcal{H} \mid \sup \scal{C-x}{y} \leq 0\}, & \text{if } x \in C; \\
   \emptyset, & \text{otherwise. } 
  \end{cases}  $$
\end{definition}

\begin{ex}\normalfont
Let $C:= \{ 0 \} \times \RR$. Let $c = (0, t) \in C$ and $ u =(u_1, u_2)$ and for all $y = (0, s) \in C$, we have $y - c = (0, s) - (0, t) = (0, s-t)$ and therefore;
\begin{align*}
\scal{u}{y - c}  & = \scal{(u_1, u_2)}{(0, s-t)} \\
& = u_1 (0) + u_2 (s-t) =  u_2 (s-t).
\end{align*}
So inorder for $u_2 (s-t) \leq 0$ for all $s \in \RR$, $u_2$ must be $0$. Hence, $$N_{C} (0,t) =  \begin{cases}
    \RR \times \{0 \}  & \text{for every } (0,t) \in C; \\
   \emptyset, & \text{otherwise. } 
  \end{cases}$$
\end{ex}
\begin{definition}\label{Def:Monotobe} \normalfont \cite[Definition 20.1 ]{BC2017} 
A set valued operator $A: \mathcal{H} \To \mathcal{H}$ is \emph{monotone} if 
$\big( \forall \big( x,  u \big) \in \gra A  \big) \big( \forall \big( y, v \big) \in \gra A  \big) $
$$   \scal{x-y}{ u - v} \geq 0. $$
\end{definition}

\begin{definition}\label{parabf}\normalfont \cite[Definition 22.1 ]{BC2017} An operator $A: \mathcal{H} \to 2^{\mathcal{H}}$ is \emph{paramonotone} if it is monotone and  $\big( \forall \big( x,  u \big) \in \gra A  \big) \big( \forall \big( y, v \big) \in \gra A  \big) $
$$   0 \leq \scal{x-y}{ u - v} = 0 \Rightarrow (x, v)  \in \gra A \ \ \text{and} \ \  (y, u)  \in \gra A.$$
\end{definition}
Paramonotone operators have attracted considerable attention in nonlinear functional analysis and optimization because of their useful properties and wide range of applications, particularly in variational inequalities, convex optimization, and equilibrium problems. These operators offer a more flexible setting than strictly monotone operators, making it possible to study solution sets and convergence properties under less restrictive assumptions. Censor in \cite{censor1998interior} highlighted their importance by establishing convergence results for an interior-point method applied to variational inequality problems involving paramonotone operators. Later, Bauschke in \cite{bauschke2012attouch} explored the connection between paramonotonicity, Attouch–Théra duality, and operator-splitting methods, showing how this property helps clarify the relationships between primal and dual solutions. In a related study, Bauschke \cite{bauschke2014rectangularity} examined the structural characteristics of paramonotone operators and their relationships with other classes of maximally monotone operators. More recently, the role of paramonotonicity has also been investigated in convex feasibility problems \cite{camacho2024feasibility}, reflecting its continued relevance in optimization research. Together, these studies show that paramonotonicity is not only an important theoretical concept but also a useful property for understanding solution structures and supporting the convergence analysis of optimization algorithms.
\begin{ex}\normalfont
Let $A = \begin{pmatrix} 0 & 0\\ 0 & 1 \end{pmatrix} $, $x = (x_1, x_2) \in \RR^{2}$ and $y   = (y_1, y_2) \in \RR^2$. Set $u_1:= x_1- y_1$ and $u_2:= x_2- y_2$. Then, 
\begin{align*}
Ax - Ay &= A(x - y)  \\
& = \begin{pmatrix} 0 & 0\\ 0 & 1 \end{pmatrix}  \begin{pmatrix} x_1-y_1\\ x_2-y_2 \end{pmatrix}  \\
& = \begin{pmatrix} 0 & 0\\ 0 & 1 \end{pmatrix}  \begin{pmatrix} u_1 \\ u_2\end{pmatrix} \\
& = \begin{pmatrix} 0 \\ u_2\end{pmatrix}.
\end{align*}
Moreover, $\scal{x- y}{Ax - Ay} = \scal{x-y}{A(x-y)} = \scal{(x_1-y_1), (x_2 - y_2)}{ (0, x_2 - y_2)} = (x_2 - y_2)^{2}  = u^{2}_2 \geq 0$. Therefore, 
\begin{equation}\label{dddo0}
\scal{x- y}{Ax - Ay}  =  u^{2}_2 \geq 0. 
\end{equation}
Furhermore, if $ \scal{x- y}{Ax - Ay}=  0$, then by \cref{dddo0}, $u^{2}_2 = 0 \Rightarrow x_2 = y_2$. Then 
\begin{align*}
Ax & = \begin{pmatrix} 0 & 0\\ 0 & 1 \end{pmatrix}  \begin{pmatrix} x_1\\ x_2 \end{pmatrix} = \begin{pmatrix} 0\\ x_2 \end{pmatrix}  = \begin{pmatrix} x_1\\ y_2 \end{pmatrix}  = Ay.
\end{align*}
Therefore, $A$ is paramontone operator in $\RR^2$. 
\end{ex}

\begin{definition}\normalfont \cite{BC2017}
A monotone operator $A: \mathcal{H} \To \mathcal{H}$ is \emph{maximally monotone} if there exists no monotone operator $B$ such that $\gra A \subset \gra B$. That is, for every $(x, u) \in \mathcal{H} \times \mathcal{H}$,
$$   \big(  \forall (x, u) \in \gra A \big) \ \big(  \forall (y, v) \in \gra A \big) \ \ \scal{x- y}{u - v} \geq 0.$$
\end{definition}
For more details about maximally monotone operators, we refer the reader to \cite{BC2017}, \cite{auslender2000asymptotic}, \cite{aborwein2010fifty}, \cite{borwein2004techniques}, \cite{burachik2002maximal}, \cite{combettes2018monotone}, \cite{eckstein1989splitting}, \cite{borwein2010convex}, and the references therein.

\begin{definition}\normalfont \cite[Definition 25.10 ]{BC2017} Let $A: \mathcal{H}\To \mathcal{H}$ be monotone. The $A$ is $3^{*}$ \emph{monotone} if 
$$  \dom A \times \ran A \subset \dom F_{A} , $$
where $F_A$ is \emph{Fitzpatrick function} of $A$.\footnote{Let $A: \mathcal{H} \to 2^{\mathcal{H}}$ be monotone. The  \emph{Fitzpatrick function} of $A$ is 
\begin{align*}
F_A: & \mathcal{H} \times \mathcal{H} \to [ - \infty, + \infty]\\
& (x, u)\mapsto \sup_{(y,v) \in \gra A} \big(   \scal{y}{u}  +\scal{x}{v} - \scal{y}{v}\big)  \\
& \scal{x}{u} - \inf_{(y,v) \in \gra A} \scal{x-y}{u-v}.
\end{align*}
 \label{footnote 1}
}
\end{definition}

\section{Parallel sums and Attouch--Th{\'e}ra Duality}
This  section reviews and records now some properties that are relevant for our studies. 
\subsection{Parallel sums} \normalfont
The \emph{parallel sum} of two operators was originally introduced by Passty \cite{passty1986parallel}. Given two set-valued operators $A$ and $B$ from $\mathcal{H}$ to $2^{\mathcal{H}}$. The \emph{parallel sum} of $A$ and $B$ is defined as,
 \begin{equation}\label{toonq}
A \Box B := \big(   A^{-1} + B^{-1}\big)^{-1},
\end{equation}
i.e.,
\begin{equation}\label{eq:001}
\big(  \forall x \in \mathcal{H} \big) \ \ \ \big(  A \Box B \big) x = \bigcup_{a+b = x}  Aa \cap  Bb.
\end{equation}

 \subsection{Recalling Attouch--Th{\'e}ra duality}
Let $A$ and $B$ be two maximally monotone operators on $\mathcal{H}$. The \emph{primal problem} associated with $ \big( A, B \big)$ is to
\begin{empheq}[box=\mybluebox]{equation}\label{Sum1}
\text{ find $x \in \mathcal{H}$ such that $0 \in Ax + Bx$.}
 \end{empheq}
 The set of \emph{primal solutions} associated with $\big( A,B \big)$ are the solutions to the corresponding sum problem \cref{Sum1} are defined as 
\begin{empheq}[box=\mybluebox]{equation}\label{SumSol}
\text{psol$( A, B):=$ zer$( A, B) = ( A, B)^{-1} (0) = \Big\{ x \in  \mathcal{H}   \mathrel{\Big|}  0 \in (A+B) x \Big\}$.}
 \end{empheq} 
 Now define 
 \begin{equation}\label{leem}
 B^{\ovee}:= (- \Id ) \circ  B \circ (- \Id)
 \end{equation}
 and $B^{-\ovee} := (B^{-1})^{\ovee} = (B^{\ovee})^{-1}$. This allows us to define the \emph{dual pair} of $(A, B)$:
 \begin{equation}
 (A, B)^{*} := ( A^{-1}, B^{- \ovee} ).
 \end{equation}
 Then the \emph{dual problem} associated with $(A, B)$ is defined to be the primal problem associated with the  dual pair $( A^{-1}, B^{-\ovee} )$:
 \begin{empheq}[box=\mybluebox]{equation}\label{Sum2}
\text{ find $y \in \mathcal{H}$ such that $0 \in A^{-1}y + B^{-\ovee}y = A^{-1} y - B^{-1} (-y)$.}
 \end{empheq} 
 The set of of \emph{dual solutions} associated with $\big( A,B \big)$ are the solutions to the corresponding sum problem \cref{Sum2}: 
\begin{empheq}[box=\mybluebox]{equation}\label{SumSol2}
\text{dsol$( A, B):=$ psol$( A, B)^{*}=$ zer $( A^{-1}+  B^{-\ovee})  = \Big\{ y \in \mathcal{H} \mathrel{\Big|}  0 \in (A^{-1}+B^{-\ovee}) y \Big\}.$}
\end{empheq}  
Because 
\begin{equation}\label{drswa11}
(A^{-1})^{-1} = A, \ \ \text{and} \ \ (A^{\ovee})^{\ovee} = A, \ \ \text{and} \ \ (A^{-\ovee})^{-\ovee} = A, 
\end{equation}
we have 
 \begin{equation}
 (A, B)^{**} = (A, B).
 \end{equation}

\section{Results}
\cref{Fact001} collects several fundamental relationships between primal and dual solution sets, operator intersections, and parallel sums. These relationships are established in the framework of Attouch–Théra duality by Bauschke, see \cite[Proposition 3.1, Corollary 3.7, and Remark 3.8]{bauschke2012attouch}. We present them in a unified formulation for use in the subsequent analysis.
\begin{lemma}\label{Fact001}\normalfont Let $A$ and $B$ be maximally monotone and assume that $x, y \in \mathcal{H}$. We recall  from \cref{SumSol} and \cref{SumSol2} that 
\begin{empheq}{equation*}
\text{psol$( A, B):=$ zer$( A, B) = ( A, B)^{-1} (0) = \Big\{ x \in  \mathcal{H}   \mathrel{\Big|}  0 \in (A+B) x \Big\}$.}
 \end{empheq} 
and 
\begin{empheq}{equation*}
\text{dsol$( A, B):=$ psol$( A, B)^{*}=$ zer $( A^{-1}+  B^{-\ovee})  = \Big\{ y \in \mathcal{H} \mathrel{\Big|}  0 \in (A^{-1}+B^{-\ovee}) y \Big\}.$}
 \end{empheq} 
Then the following holds;
\begin{enumerate}
\item\label{gs1} $ \dom (A \cap (-B) )= \text{psol} ( A, B)$.
\item\label{gs2} $\ran (A \cap (-B) ) = \text{dsol} ( A, B)$.
\item\label{gs3} $ \dom (A^{-1} \cap B^{-1} (- \Id) )= \text{dsol} ( A, B)$
\item\label{gs4} $ \ran (A^{-1} \cap B^{-1} (- \Id) )= \text{psol} ( A, B)$
\item\label{gs5} $\gra (A \cap (-B) )$ and $\gra (A^{-1} \cap B^{-1} (- \Id) ) $ are closed and convex sets.
\item\label{gs6} $ \text{psol} ( A, B) = ( A+B)^{-1} (0)= (A^{-1} \Box B^{-1})(0).$
\item\label{gs7} $ \text{dsol} ( A, B) = ( A^{-1}+B^{-\ovee})^{-1} (0)= (A \Box B^{-\ovee})(0).$
\end{enumerate}
\end{lemma}
\begin{proof}
\ref{gs1}: "$\subseteq$": Let  $\emptyset \neq \dom (A \cap (-B) )$. Then $\exists \ y$ such that $y \in Ax$ and  $y \in - B x$. This yelids to $y \in Ax $ and  $- y \in  B x$. Therefore, 
$0 = y+ (-y)  \in Ax +B x$. Hence,  $0  \in (A +B) x$. This implies that $ x  \in (A +B)^{-1} (0)$. Therefore,  $x  \in \text{psol} ( A, B)$. This gives that 
\begin{equation}\label{ans1}
\dom (A \cap (-B) ) \subseteq \text{psol} ( A, B).
\end{equation}
Next, "$\supseteq$": Take $x \in \text{psol} ( A, B)$. This implies that $ x \in (A+B)^{-1} (0)$. Therefore, $0 \in (A+B) x$. Then there exists $a, b \in \mathcal{H}$ such that $a \in Ax$ and $b \in B x$ such that $a + b = 0$. This implies that $a \in Ax$ and $-a \in Bx$, where $b = -a$. Hence, $a \in Ax$ and $a \in -Bx$. Sequencely, $a \in (A \cap -B)(x)$. 
This implies that $ \dom (A \cap -B) \neq \emptyset$. Therefore, we obtain
\begin{equation}\label{ans2}
\dom (A \cap (-B) )\supseteq \text{psol} ( A, B).
\end{equation}
From \cref{ans1} and \cref{ans2} we have 
$$ \dom (A \cap (-B) ) = \text{psol} ( A, B).   $$
\ref{gs2}: Note that $ \ran (A \cap (-B) ) = \{ y \in \mathcal{H}: \exists x \in \mathcal{H}  \ \ \text{such that} \ \ y \in (A \cap (-B)) (x)  \}$. So let $y \in \ran (A \cap (-B)) \Leftrightarrow \exists \ x \in \mathcal{H}$ such that $y \in Ax \cap (-Bx)$ .This is equivalent to $ \exists y \in \mathcal{H}$ such that $y \in Ax$ and $y \in ( (- \Id ) \circ  B \circ (- \Id) (-x))$ by \cref{leem}. This yelids to $\exists y \in \mathcal{H}$ such that $y \in Ax$ and $y \in B^{\ovee} (-x) \Leftrightarrow \exists y \in \mathcal{H} \ \ \text{such that} \ \ x \in A^{-1} y \ \ \text{and} \ \ -x \in B^{- \ovee} y \Leftrightarrow \exists y \in \mathcal{H} \ \ \text{such that} \ \ 0 = x + (-x) \in ( A^{-1} + B^{- \ovee} ) y.$ This gives, $y \in \text{dsol} ( A, B) $. Hence 
$$ \ran (A \cap (-B) ) = \text{dsol} ( A, B).   $$
\ref{gs3} and \ref{gs4}: The proofs are identical to that of  \ref{Fact01}\ref{gs1} and 
 \ref{Fact01}\ref{gs2}, respectively. 
\ref{gs5}: See \cite[Proposition 20.21 \& 20.33]{BC2017} as well as \cite[Theorem~4.6]{klee1970convex}.
\ref{gs6}: By using \cref{toonq} and \cref{drswa11}, we have 
\begin{align*}
(A^{-1} \Box B^{-1})(0) & = \big( (A^{-1})^{-1} + ( B^{-1})^{-1} \big)^{-1} (0)\\
& =  \big( A + B  \big)^{-1} (0) \\
& = \text{psol} ( A, B) \ \ \ \text{by \cref{SumSol}}.
\end{align*}
\ref{gs7}: By using \cref{toonq},  \cref{drswa11}, and \cref{SumSol2} the proof is the same to that of \cref{gs6}.
\end{proof}

\begin{lemma}\label{lojec1}\normalfont 
Let $A, B: \mathcal{H} \To \mathcal{H}$ and let $(x, u) \in \mathcal{H} \times \mathcal{H}$. Then 
 $(x, u) \in \gra (A^{-1} + B^{-1})^{-1} \Leftrightarrow (u,x) \in \gra (A^{-1} + B^{-1}). $
\end{lemma}
\begin{proof}
Let $(x, u) \in  \gra (A^{-1} + B^{-1})^{-1} \Leftrightarrow u \in (A^{-1} + B^{-1})^{-1} x \Leftrightarrow x \in (A^{-1} + B^{-1}) u.$ Therefore, $(u, x) \in \gra (A^{-1} + B^{-1})$.
\end{proof}

\begin{lemma}\label{lojec22}
Let $A, B: \mathcal{H} \To \mathcal{H}$ and let $( \bar{x}, \bar{u}) \in \mathcal{H} \times \mathcal{H}$. Then 
$$ ( \bar{x}, \bar{u}) \in \gra (A^{-1}+ B^{-1}) \Leftrightarrow (u,  x)  \in \gra (A^{-1}+ B^{-1})^{-1} - \Id, \ \text{where} \ \ u = \bar{u}- \bar{x} \ \text{and} \  x=  \bar{x}.$$
\end{lemma}
\begin{proof}
Let $ (x, u) \in \gra (A^{-1} + B^{-1})^{-1}- \Id$. Then \cref{lojec1} gives $(x, u) \in \gra (A^{-1} + B^{-1})^{-1}- \Id \Leftrightarrow u+ x \in  (A^{-1} + B^{-1})^{-1} \Leftrightarrow x \in (A^{-1} + B^{-1}) (u+ x) \Leftrightarrow \bar{x} \in (A^{-1} + B^{-1}) ( ( \bar{u} - \bar{x}) + \bar{x}) \Leftrightarrow  \bar{x} \in (A^{-1} + B^{-1}) (\bar{u} ) \Leftrightarrow  ( \bar{x}, \bar{u}) \in \gra (A^{-1}+ B^{-1}). $
\end{proof}

\begin{theorem}\label{Fact01}\normalfont  Let $A$ and $B$ be set-valued operators on $\mathcal{H}$. Then the following hold:
\begin{enumerate}
\item\label{01} $ \dom \big(  A \Box B \big) = \ran \big(  A^{-1}+ B^{-1}\big) \subseteq \ran A^{-1} + \ran B^{-1} = \dom A+ \dom B.$
\item\label{02} $ \ran \big(  A \Box B \big) = \ran A \cap \ran B$.
\item\label{03} If $A$ and $B$ are monotone, then so is $A \Box B$.
\item\label{003} Suppose that $A^{-1} + B^{-1}$ is maximally monotone and define $E:= (A^{-1} + B^{-1} )^{-1}- \Id. $ Then $E$ is \emph{maximally} $1$-\emph{hypomonotone}, in the sense that $E+ \Id$ is maximally monotone. 
\item\label{04} If $A^{-1} + B^{-1}$ is maximally monotone, then so is $A \Box B$.
\end{enumerate}
\end{theorem}

\begin{proof}
\ref{01}: By using \cref{toonq}, let $x \in  \dom \big(  A \Box B \big) \Leftrightarrow \ x \in  \dom \big(  A^{-1} + B^{-1} \big)^{-1} \Leftrightarrow \exists \ y \in \mathcal{H}$ such that $y \in (A^{-1} + B^{-1} )^{-1} x \Leftrightarrow x \in (A^{-1} + B^{-1} ) y \Leftrightarrow x \in \ran (A^{-1} + B^{-1} )$. Moreover, let $y \in \ran (A^{-1} + B^{-1} ) (x),$ then, $y \in A^{-1}(x)+ B^{-1} (x)$. Hence, $\exists a \in A^{-1} (x) \Rightarrow a \in \ran (A^{-1})$ and $b \in B^{-1} (x) \Rightarrow b \in \ran (B^{-1}) $ such that $y = a+ b $. Therefore, 
$$  y = a + b \in \ran (A^{-1}) + \ran (B^{-1}).$$
Since $\ran A^{-1} = \dom A$ and  $\ran B^{-1} = \dom B$, then by using Minkowski sum we obtain $\ran A^{-1} + \ran B^{-1} = \dom A+ \dom B.$
\ref{02}: "$\Rightarrow$ Let $y \in \ran (A \Box B ) \Leftrightarrow \exists x \in \mathcal{H}$ such that $ y \in (A \Box B ) x \Leftrightarrow \exists a, b \in \mathcal{H}$ such that $ x = a+ b$ and $a \in A^{-1} y \Leftrightarrow y \in A(a) \Leftrightarrow y \in \ran (A)$ and $b \in B^{-1} y \Leftrightarrow y \in B(b) \Leftrightarrow y \in \ran (B)$. Therefore, $ y \in \ran A \cap  \ran B \Leftrightarrow ran (A \Box B) \subseteq \ran A \cap \ran B$. "$\Leftarrow$" Let $ y \in \ran A \cap  \ran B \Leftrightarrow y \in \ran A$ and $ y \in \ran B \Leftrightarrow A^{-1} (y) \neq \emptyset$ and $B^{-1} (y) \neq \emptyset$. Hence, $ \exists \ a, b \in \mathcal{H}$ such that $a \in A^{-1} (y)$, $b \in B^{-1} (y)$, and
$ a+ b \in ( A^{-1}  + B^{-1} ) (y) \Leftrightarrow y \in \dom (A^{-1} + B^{-1})  \Leftrightarrow y \in  \ran (A^{-1} + B^{-1})^{-1} \Leftrightarrow y \in A \Box B$, by \cref{toonq}.
\ref{03}: Take $( x_1, y_1) $ and $( x_2, y_2) $ in $ \gra \big(  A^{-1}+ B^{-1} \big)^{-1} $. This implies that $y_1 \in \big(  A^{-1}+ B^{-1} \big)^{-1} x_1 \Leftrightarrow x_1 \in  A^{-1}y_1 + B^{-1} y_1 $ and $y_2 \in \big(  A^{-1}+ B^{-1} \big)^{-1} x_2 \Leftrightarrow x_2 \in A^{-1} y_2 + B^{-1} y_2$. Then there exists $ a_1, a_2, b_1, b_2 \in \mathcal{H}$ such that $x_1 = a_1 + b_1$ and $x_2 = a_2 + b_2$, where $a_1 \in A^{-1} y_1, b_1 \in B^{-1} y_1, a_2 \in A^{-1} y_2,$ and $b_2 \in B^{-1} y_2$. By using \cref{Def:Monotobe}, we obtain 
\begin{align*}
\scal{x_1- x_2}{y_1 - y_2} & = \scal{(a_1+b_1) - (a_2+b_2)}{y_1 - y_2}  \\
& =   \scal{(a_1-a_2) + (b_1-b_2)}{y_1 - y_2}  \\
& =   \scal{a_1-a_2}{y_1 - y_2}  + \scal{b_1-b_2}{y_1 - y_2} \\
&  \geq 0+0 =0.
\end{align*}
\ref{003}: Let $(x, u) \in \mathcal{H} \times \mathcal{H}$  satisfies 
$$ \forall (y, v) \in \gra E  \ \ \  \scal{x-y}{u-v}  \geq - \norm{x - y}^2 $$
This implies that 
\begin{align*}
& \scal{x-y}{u-v}  \geq - \norm{x - y}^2 \\
& \scal{x-y}{u-v} + \scal{x-y}{x-y} \geq 0 \\
& \scal{x-y}{ (u +x) -(v + y)}  \geq 0 
\end{align*}
By using \cref{lojec22}  and the definition of $E$, we obtain $ (y, v) \in \gra E \Leftrightarrow (y, v+ y) \in \gra (A^{-1} +B^{-1} )^{-1}$. Therefore, $\forall (y, c) \in \gra (A^{-1} +B^{-1} )^{-1} \ \  \scal{x-y}{ (u +x) - c }  \geq 0$, where $c:= v+y$. Note that, $(A^{-1} +B^{-1} )^{-1}$ is maximally monotone. This gives, $(x, u+ x ) \in \gra (A^{-1} +B^{-1} )^{-1} \Leftrightarrow ( u + x) \in (A^{-1} +B^{-1} )^{-1} x  \Leftrightarrow u \in (A^{-1} +B^{-1} )^{-1} x - x \Leftrightarrow u \in ( (A^{-1} +B^{-1} )^{-1} - \Id ) x \Leftrightarrow u \in E (x)$. Therefore, $ (x , u) \in \gra E$.
\ref{04}: Since $\Id$ is maximally motone and full domain and  $(A^{-1} + B^{-1})^{-1} - \Id $ is maximally $1$-hypomonotone by \ref{003}. Hence,  by \cref{toonq} we have $(A^{-1} + B^{-1})^{-1} - \Id  + \Id = (A^{-1} + B^{-1})^{-1} = A \Box B,$ is maximally motone. 
\end{proof}

\begin{theorem}\label{Fact02}\normalfont  Let $A$ and $B$ be set-valued operators on $\mathcal{H}$. Then the following hold:
\begin{enumerate}
\item\label{11} If $A$ and $B$ are $3^{*}$ monotone, then so is $A \Box B$.
\item\label{12} If $A$ and $B$ are paramonotone, then so are $A^{-1}, B^{-1}$ and $ A \Box B$.
\item\label{13} If $A$ and $B$ are maximally monotone, then for every $x \in X$, the set $ \big( A \Box B \big) x$ is convex.
\item\label{14} If $ 0 \in \sri \big(  \dom f^* - dom g^*\big)$, then $ \big( \partial f \big) \Box \big( \partial g \big) = \partial \big( f \Box g \big)$ and the infimal convolution is exact everywhere.
\end{enumerate}
\end{theorem}

\begin{proof}
\ref{11}: It follows directry from  \cite[Proposition 25.22]{BC2017}.  \ref{12}: See \cite[Proposition 22.2]{BC2017} for more details. 
\ref{13}: It comes from \cite[Corollary 3.4]{bauschke2012attouch}
\ref{14}: See \cite[Proposition 25.23]{BC2017} .
\end{proof}

The next \cref{Lem1} provides a direct formulation of this invariance principle in terms of the parallel sum, without requiring maximal monotonicity for the principal equality. This formulation may be useful in subsequent investigations involving parallel sums, normal cone operators, and related constructions.
\begin{lemma}\label{Lem1}\normalfont Let $A$ and $B$ be paramonotone on $X$, and let $x \in X$. Suppose that $a_1, a_2, b_1, b_2$ are vectors in $X$ such that $a_1+ b_1 = x = a_2 + b_2$ and $A a_1 \cap Bb_1 \neq  \emptyset \neq Aa_2 \cap Bb_2$. Then 
\begin{equation}\label{eq:003}
\big(  A \Box B \big) x = A a_1 \cap B b_1 = A a_2 \cap B b_2;
\end{equation}
consequently, if $A$ and $B$ are also maximally monotone, then 
\begin{equation*}
\big(  A \Box B \big) x \ \ \text{is closed and convex.}
\end{equation*}
\end{lemma}

\begin{proof}
Let $u_1 \in A a_1 \cap B b1 \Rightarrow u_1 \in A a_1$ and $u_1 \in B b_1$. Also, $u_2 \in A a_2 \cap B b2 \Rightarrow u_2 \in A a_2 $ and $u_2 \in B b_2 $. Since $A$ and $B$ are monotone, then 

\begin{equation}\label{eq1}
\big( \forall \big(  a_1, u_1\big) \in \gra A  \big) \big( \forall \big(  a_2, u_2\big) \in \gra A  \big)  \ \ \ \ 0 \leq \scal{a_1-a_2}{ u_1-u_2}, 
\end{equation}
and 
\begin{equation}\label{eq2}
\big( \forall \big(  b_1, u_1\big) \in \gra B  \big) \big( \forall \big(  b_2, u_2\big) \in \gra B  \big)  \ \  \ \ 0 \leq \scal{b_1-b_2}{ u_1-u_2}.
\end{equation}
Adding \cref{eq1} and \cref{eq2} gives
\begin{align*} 
0 \leq \scal{a_1-a_2}{ u_1-u_2} + \scal{b_1-b_2}{ u_1-u_2}& = 2 \scal{ \big( a_1 + b_1 \big) -\big( a_2 + b_2\big) }{ u_1-u_2} \\
& = 2 \scal{x-x}{ u_1-u_2}   \ \ (\text{from assumption})\\
& = 0.
\end{align*}
Hence, 
$$ \scal{a_1-a_2}{ u_1-u_2}  = 0 \ \ \ \text{and}  \ \ \  \scal{b_1-b_2}{ u_1-u_2}  = 0.$$
But $A$ is paramonotone so $$u_2 \in Aa_1  \ \ \ \text{and} \ \ \  u_1 \in A a_2. $$
Similary, because $B$ is paramonotone, we obtain 
$$u_2 \in Bb_1  \ \ \ \text{and} \ \ \  u_1 \in B b_2. $$
Altogether, 
$$ u_2 \in Aa_1 \cap Bb_1 \ \ \ \text{and} \ \ \  u_1 \in A a_2 \cap Bb_2 .$$
Therefore, we have 
$$  Aa_1 \cap Bb_1  = A a_2 \cap Bb_2 .$$
This is true for all $x$, then by \cref{eq:001} we have 
$$ \big( A \Box B \big) x =  Aa_1 \cap Bb_1  = A a_2 \cap Bb_2 , $$
which verifies \cref{eq:003}.
\end{proof}

The next theorem examines the parallel sum of a set-valued operator $A$ and the normal cone operator $N_{V}$ associated with a closed linear subspace $V$. We first derive explicit characterizations of its values, domain, and range using the standard properties of parallel sums and normal cone operators. We then establish, by a direct argument based on paramonotonicity and orthogonality, that the value of the parallel sum can be determined from any admissible decomposition. This provides a specialized formulation of the decomposition-invariance phenomenon associated with paramonotone operators, closely related to the results of Bauschke, see \cite{bauschke2012attouch} for more details.

\begin{theorem}\label{theorem111}\normalfont 
Let $A$ be a set-valued operator on $X$ and let $V$ be a closed linear subspace of $X$. Then the following hold; 
\begin{enumerate}
\item\label{gg1s} $ \big(  \forall x \in X \big) \  \ \big(  A \Box N_{V} \big) \big( x\big) = V^\perp \cap A \big( x - V \big)$. 
\item\label{gg2s} $ \dom \big(  A \Box N_{V}\big) = V + A^{-1} \big( V^{\perp}\big).$
\item\label{gg233} $ \ran \big(  A \Box N_{V}\big) = \ran A \cap V^{\perp}$.
\end{enumerate}
Additionally, if  $A$ is paramonotone and $x \in \dom \big(  A \Box N_{V}\big)$, i.e., $x = a+ v$ where $Aa \cap V^{\perp} \neq \emptyset$ and $v \in V$, then
\begin{equation}\label{gfgdatrrrr}
\big(  A \Box N_{V} \big) \big( x \big) = Aa \cap V^{\perp}.
\end{equation}
\end{theorem}

\begin{proof}
Let $x \in \mathcal{H}$ and $w \in \big(  A \Box N_{V}  \big) \big( x\big) $. This equivalence to
\begin{align*}
 \ \ \ \ \ w & \in \big(  A \Box N_{V}  \big) \big( x\big) \\
\Leftrightarrow w & \in Aa \cap N_{V} \big( v \big) \ \ \text{for $x = a+ v$ and $v \in V$  and  by \eqref{..}} \\
\Leftrightarrow w  &  \in Aa \cap V^\perp \\
\Leftrightarrow w  &  \in A (x - V) \cap V^\perp \ \ \text{for all } \ \ v \in V \ \  \text{and} \ \ a = x-v.
\end{align*}
Next, let $x \in  \dom \big(  A \Box N_{V}\big) $. This is equivalent to 
\begin{align*}
 \ \ \ \ \  & x  \in  \dom \big(  A \Box N_{V}\big)  \\
\Leftrightarrow x & \in  \dom \big(  A (x - V) \cap V^\perp \big)  \ \ \  by \ \  \text{\cref{theorem111}\ref{gg1s} } \\
\Leftrightarrow  \ \ \ & \dom \big(  A (x - V) \cap V^\perp \big) \neq \emptyset \ \ \  \\
\Leftrightarrow v^{\perp}  &  \in A (x - V) \cap V^\perp \ \ \text{for some } \ \ v \in V \ \  \text{and} \ \ v^{\perp} \in V^{\perp}\\
\Leftrightarrow v^\perp & \in  A(x- V) \ \  \text{for some } \ \ v \in V \ \  \text{and} \ \ v^{\perp} \in V^{\perp} \\
\Leftrightarrow x  &  \in v+ A^{-1} (v^{\perp}) \ \  \text{for some } \ \ v \in V \ \  \text{and} \ \ v^{\perp} \in V^{\perp} \\
\Leftrightarrow x  &  \in v+ A^{-1} (V^{\perp}) . 
\end{align*}
Hence, 
$$ \dom \big(  A \Box N_{V}\big) = V + A^{-1} \big( V^{\perp}\big).  $$
Next, assume that $A$ is paramonotone and $x \in \dom \big(  A \Box N_{V}\big)$, say, $x = a+ v$ where $Aa \cap V^{\perp} \neq \emptyset$ and $v \in V$, then applying \cref{Lem1} gives that 
\begin{align*}
\big(  A \Box N_{V} \big) \big( x \big) & = A(a)  \cap N_{V}(v) \\
& =  A(a)  \cap V^{\perp}.
\end{align*}
\end{proof}

\begin{remark}\normalfont 
The identities in \cref{theorem111}\ref{gg1s}-\ref{gg233}  follow from the definition of the parallel sum and the structure of the normal cone operator of a closed linear subspace. The additional representation in \eqref{gfgdatrrrr} follows from paramonotonicity and the orthogonality relation between $V$ and $V^\perp$. It may also be viewed as a specialization of the general decomposition-invariance property for paramonotone parallel sums. Our presentation emphasizes a direct proof under the assumption of paramonotonicity, without requiring maximal monotonicity of $A$.
\end{remark}

\section*{Acknowledgments}
The author thank the Royal Commission for Jubail and Yanbu for supporting her reaseach.

\end{document}